\pdfoutput=1
\documentclass[letterpaper,draftcls,onecolumn]{ieeeconf} 
\usepackage{graphicx}
\usepackage{amsmath}
\usepackage{amsfonts}
\usepackage[psamsfonts]{amssymb}
\usepackage{times}
\usepackage{cases}
\usepackage{color}
\usepackage{array}
\usepackage{mathrsfs}
\usepackage{mathtools}
\let\labelindent\relax
\usepackage{enumitem}

\mathtoolsset{showonlyrefs}

\newtheorem{definition}{Definition}

\newtheorem{lemma}{Lemma}

\newtheorem{theorem}{Theorem}
\newtheorem{talgorithm}{Algorithm}
\newtheorem{assumption}{Assumption}

\newtheorem{problem}{Problem}

\makeatletter
\newcommand{\pushright}[1]{\ifmeasuring@#1\else\omit\hfill$\displaystyle#1$\fi\ignorespaces}
\newcommand{\pushleft}[1]{\ifmeasuring@#1\else\omit$\displaystyle#1$\hfill\fi\ignorespaces}
\makeatother

\newcommand{\R}{\mathbb{R}}
\newcommand{\N}{\mathbb{N}}

\newcommand{\p}{\mathbb{P}}

\newcommand{\barx}{\bar{x}} 
\newcommand{\xopt}{\hat{x}} 
\newcommand{\mopt}{\hat{\mu}} 
\newcommand{\uinp}{\tilde{\ell}^{s_i}_{p_i}}
\newcommand{\lspi}{\ell^{s_i}_{p_i}}
\newcommand{\inp}{\tilde{\ell}^s_p}
\newcommand{\trunc}{P_{t}}
\newcommand{\adj}{\textnormal{Adj}}

\newcommand{\sys}{\mathcal{G}}

\newcommand{\outp}{\tilde{\ell}^r_q}
\newcommand{\Mjdp}{\mathcal{M}}

\newcommand{\bx}{x}
\newcommand{\btx}{\tilde{x}}

\newcommand{\lap}{\textnormal{Lap}}
\newcommand{\Lap}{\textnormal{Lap}}
\newcommand{\runtime}{k}

\newcommand{\sone}{x^1_i}
\newcommand{\stwo}{x^2_i}

\newcommand{\fgaini}{\rho_i}
\newcommand{\lipfi}{K_{i}}

\newcommand{\lipg}{K_g}

\newcommand{\lipgradgi}{L_{g,i}}

\newcommand{\E}{\mathbb{E}}
\newcommand{\tx}{\tilde{x}}

\newcommand{\M}{M}

\author{M.T. Hale and M. Egerstedt$^\star$\thanks{$^\star$The authors are with the School of Electrical and Computer Engineering, Georgia Institute of
Technology, Atlanta, GA 30332, USA. Email: \texttt{\{matthale, magnus\}@gatech.edu}. Research supported in
part by  the NSF under Grant CNS-1239225.}
}

\title{Approximately Truthful Multi-Agent Optimization Using Cloud-Enforced Joint Differential Privacy}

\begin{document}
\maketitle
\thispagestyle{empty}
\pagestyle{empty}

\begin{abstract}
Multi-agent coordination problems often require 
agents to exchange state information in order to reach
some collective goal, such as agreement on a final
state value. In some cases, it is feasible that opportunistic
agents may deceptively report false state values for their own benefit,
e.g., to claim a larger portion of shared resources.
Motivated by such cases, this paper presents a multi-agent
coordination framework which disincentivizes opportunistic misreporting of
state information. This paper focuses on multi-agent coordination problems that
can be stated as nonlinear programs, with non-separable constraints coupling
the agents. In this setting, an opportunistic agent may be tempted to 
skew the problem's constraints in its favor to reduce its 
local cost, and this is exactly the behavior we seek to disincentivize. 
The framework presented uses a primal-dual approach
wherein the agents compute primal updates and a centralized cloud
computer computes dual updates. 
All computations performed by the cloud are carried out in a way that
enforces \emph{joint differential privacy}, which adds noise in order
to dilute any agent's influence upon the value of its cost function in the
problem. We show that this dilution deters agents from 
intentionally misreporting their states to the cloud, and present bounds
on the possible cost reduction an agent can attain through misreporting its state.
This work extends our earlier work on incorporating ordinary 
differential privacy into multi-agent optimization, and we show
that this work can be modified to provide a disincentivize for misreporting states to the cloud.
Numerical results are presented to demonstrate convergence of the optimization
algorithm under joint differential privacy. 
\end{abstract}

\section{Introduction}
Multi-agent optimization problems have found applications
in a range of research areas, including power systems 
in \cite{nazari14},
machine learning in \cite{boyd11}, sensor networks in \cite{cortes02}, 
and robotics in \cite{soltero13}. 
Solutions to some problems in these applications rely on the implicit assumption
that all agents share correct, truthful information with others in a network. However,
one can envision a scenario in which deceitful agents do not honestly share their data,
instead sharing false information that will skew the behavior of the system in their
favor. 
For example, a homeowner connected to a smart power grid may report a false value for his or
her power usage in order to save money. 
This paper considers optimization problems with 
agents that may be intentionally deceitful
for their own benefit, and it provides a method for disincentivizing untruthful behavior
when teams of agents are collectively optimizing. 

We reduce the incentive to share false information by using
\emph{joint differential privacy}, defined in \cite{kearns12}, to limit the possible decrease in cost an agent
can achieve through intentionally misreporting its state. Joint differential privacy adds noise to reduce
the ability of an agent to benefit from reporting false information, and 
was first introduced for this purpose in \cite{kearns12} to promote
truthful sharing of information by the players in a class of games.
This work 
was applied specifically to optimization problems arising from 
 distributed electric vehicle charging in \cite{han15}, as well as to
linearly separable optimization problems in \cite{hsu16}. 

The developments in this paper differ from those in \cite{han15} and \cite{hsu16}, as well
as some other work on private optimization, e.g., \cite{han16}, in two key ways. 
First, the work in \cite{han15}, \cite{hsu16}, and \cite{han16} uses differential privacy as it was
originally defined in \cite{dwork06c} to keep some static object, such as a constraint function, private. 
Our work focuses on trajectory-level privacy for state trajectories
that are dynamically generated as an optimization algorithm is executed. 
In applications such as phase 
synchronization in smart power grids, each iterate of an optimization
algorithm corresponds to a physical state at some point in time, and an agent's
contribution to the optimization process is its whole state trajectory. 
Applying joint differential privacy in such applications should therefore be done
at the trajectory level, and implementing joint differential privacy for trajectory-level
data is most naturally done using the dynamical systems formulation of differential privacy.
Therefore, this paper uses trajectory-level privacy as defined in \cite{leny14}, rather than
privacy for databases as in \cite{dwork06c}. 

This form of privacy is not only better suited to the applications of interest, 
but also allows privacy guarantees to hold across infinite time horizons which, in many cases, cannot
be attained using privacy for databases. 
To our knowledge, this paper is the first use of joint differential privacy 
at the trajectory level. The second key difference between the existing literature and
this paper is the class of optimization problems solved. 
We apply joint differential privacy
in a general multi-agent nonlinear programming setting,
incorporating non-separable functional constraints that are also possibly nonlinear, which
differs from work done in \cite{han15}, \cite{hsu16}, and \cite{han16} where either linear or affine constraints are considered. 

To solve such problems, a cloud computer is added to the team of agents to serve as 
a trusted central aggregator, also sometimes called a ``curator'' in the privacy literature. 
This architecture was previously used for privately solving multi-agent
nonlinear programs in \cite{hale15c},  
in which honest agents seek to protect
sensitive data from eavesdroppers via ordinary differential privacy. The work
in \cite{mcsherry07} suggests that any
differential privacy implementation provides some disincentive against
untruthful information sharing, and the current paper uses the notion of joint differential privacy
to formally provide this guarantee in the framework developed in \cite{hale15c}.
Accordingly, the technical novelty of this paper is not in the algorithm used to solve problems,
but in the theoretical performance guarantees that are provided using this framework; 
the work in \cite{hale15c} focuses exclusively on protecting sensitive data, 
but the current paper focuses on the problem of preventing untruthful behavior by the agents. 
The technical contribution of this paper thus consists of adapting our existing privacy framework to
the problem of incentivizing truthful behavior by the agents, and quantifying the extent to which
any agent can benefit from untruthful behavior. 


Several existing approaches use behavioral analysis to 
identify untruthful agents, including those in \cite{braynov04} and \cite{fagiolini08}. In this paper, an agent's
local state updates rely on a local objective function and local constraint set that
are considered sensitive, and therefore these local data are not shared with any other agent. 
As a result, 
the correct next value of an agent's state is only known to that agent, and
a behavioral analysis approach cannot be used here because no outside observer can
determine what any agent's future states should be. 
Rather than detecting untruthful behavior, this paper seeks
to prevent it outright by reducing the incentive to share untruthful information via
joint differential privacy. 
The principle underlying this application of joint differential privacy is that adding noise
to the system can make each agent's cost insensitive to changes in the agent's state trajectory.  
The amount of noise added must be calibrated to the system
in order to dilute the effect of an agent reporting an untruthful
state to the cloud, and this paper presents this calibration for
a general multi-agent nonlinear program in terms of constants pertaining to the problem. 

The problems considered consist of a collection of agents, each with a local
objective function and local set constraint, 
and ensemble state constraints that jointly
constrain the agents. A primal-dual approach is used, in which
the agents update their own states, which are the problem's primal
variables, and the cloud updates the problem's dual variables. 
Naturally, the constraints in this problem will usually lead to higher costs for each
agent than a comparable unconstrained problem.
As a result, some of the agents may wish to skew the constraints in their favor.
One way they may do so is by reporting
false state information to the cloud
in order 
to loosen the constraints' effects on their on their own states,
thereby giving the untruthful agents a lower cost. 
This manipulation of the constraints will result in an unequal distribution of 
the burden of these constraints by tightening them on honest agents. Therefore,
this form of untruthful behavior is disincentivized using 
joint differential privacy. 

The rest of the paper is organized as follows. Section II provides the 
necessary optimization background, the structure of communications
in the system, and formally states the joint differentially private
optimization problem that is the focus of the paper. Then, Section III reviews 
joint differential privacy, and Section IV presents the proposed joint differentially private algorithm. 
Next, Section V proves the main result of the paper on limiting an agent's incentive to misreport
its states to the cloud.  
Section VI then presents simulation results
and Section VII concludes the paper.

\section{Background on Optimization and Problem Statement} \label{sec:problem}
This section presents the multi-agent optimization 
problem of interest and an optimization algorithm
that will later be used to solve a joint differentially private
version of this problem. 
This section also describes the cloud-based architecture
used to solve this problem. Throughout the paper, 
$\nabla_{x_i}$ denotes the partial derivative
with respect to $x_i$, and $h_{x_i}$ denotes
the partial derivative of the function $h$ with respect to $x_i$. 

\subsection{Optimization Problem Formulation}
Consider a problem consisting of $N$ agents indexed over
$i \in [N] := \{1, \ldots, N\}$. 
Agent $i$
has state $x_i \in \R^{n_i}$ with $n_i \in \N$ and a local
set constraint of the form
$x_i \in X_i \subset \R^{n_i}$.
The diameter of $X_i$ is denoted by $D_i := \max_{x_1, x_2 \in X_i} \|x_1 - x_2\|_1 $. 
Agent $i$ also has a local objective function $f_i : \R^{n_i} \to \R$ 
depending only upon its own state. Agent $i$'s local data are subject to the following assumption. 

\begin{assumption} \label{as:fX}
For all $i \in [N]$, $f_i$ is $C^2$ and convex in $x_i$, and $X_i$ is non-empty, compact, and convex. \hfill $\triangle$
\end{assumption}

Assumption \ref{as:fX} implies that $f_i$ is Lipschitz and its Lipschitz\footnote{All
Lipschitz constants in this paper are with respect to the metric induced by the $1$-norm.}
constant is denoted
by $K_i$. 
In particular, Assumption \ref{as:fX} allows for convex polynomial objectives and
box constraints, which are common
in multi-agent optimization problems. Both $f_i$ and $X_i$ are considered sensitive information
and are therefore not shared with the other agents or with the cloud. 
For simplicity of notation, define the set
$X := X_1 \times \cdots \times X_N \subset \R^n$,
where $n = \sum_{i \in [N]} n_i$.
The agents' individual set constraints require $x \in X$, 
where $x$ is the ensemble state vector of the network, defined as
$x = \left(x_1^T,  \ldots x_N^T\right)^T \in X$,
and where Assumption \ref{as:fX} provides that $X$ is non-empty, compact, and convex. 

The agents are together subject to global inequality constraints $g(x) \leq 0$,
where $g : \R^n \to \R^m$. 
The functional constraints in $g$ are subject to the following assumption. 
\begin{assumption} \label{as:g} 
\leavevmode
\begin{enumerate}[label={\upshape\roman*.}, align=left, widest=ii]
\item For all $j \in [m]$, the constraint $g_j$ is $C^2$ and convex in $x$.\label{as:gi}
\item There exists a point $\barx \in X$ such that $g(\barx) < 0$\label{as:gii}. \hfill $\triangle$
\end{enumerate}
\end{assumption}

Assumption \ref{as:g}.\ref{as:gi} admits a wide
variety of constraint functions, e.g., any convex polynomials.
Assumption \ref{as:g}.\ref{as:gii} 
is known as Slater's condition, e.g., Assumption 6.4.2 in \cite{bertsekas03},
and, in conjunction with Assumption \ref{as:g}.\ref{as:gi}, guarantees 
that strong duality holds. 
Assumption \ref{as:g} implies that $g$ and each $g_{x_i}$ are Lipschitz.
Their Lipschitz constants are denoted by $\lipg$ and $\lipgradgi$, respectively. 

Summing the per-agent objective functions gives the ensemble objective
$f(x) = \sum_{i \in [N]} f_i(x_i)$,
which is $C^2$ and convex in $x$ because of Assumption \ref{as:fX}. 
Together $f$, $g$, and $X$ comprise the following ensemble-level optimization problem.
\addtocounter{problem}{-1}
\begin{problem} \label{prob:ensemble} \emph{(Preliminary; no joint differential privacy yet)}
\begin{align*}
\textnormal{minimize }  &f(x) \\
\textnormal{subject to } &g(x) \leq 0 \\
                        &x \in X. \tag*{$\lozenge$}
\end{align*}
\end{problem} 

We now detail how to solve Problem \ref{prob:ensemble} in the cloud-based system, and
then give a unified problem statement, including incentivizing truthful behavior, in
Problem~\ref{prob:full} below. 

The Lagrangian associated with Problem \ref{prob:ensemble} is
\begin{equation}
L(x, \mu) = f(x) + \mu^Tg(x),
\end{equation}
where $\mu \in \R^m_{+}$ is the dual vector associated with Problem~\ref{prob:ensemble}
and $\R^m_{+}$ denotes the non-negative orthant of $\R^m.$
Seminal work of \cite{kuhn51}
shows that
a point $\xopt \in X$ is a solution
to Problem \ref{prob:ensemble} if and only if there exists a point $\mopt \in \R^m_{+}$
such that $(\xopt, \mopt)$ is a saddle point of $L$. This saddle
point condition can be compactly expressed by requiring
\begin{equation} \label{eq:saddle}
L(\xopt, \mu) \leq L(\xopt, \mopt) \leq L(x, \mopt) \textnormal{ for all } (x, \mu) \in X \times \R^m_{+},
\end{equation}
and an optimal primal-dual pair $(\hat{x}, \hat{\mu})$ exists under Assumptions \ref{as:fX} and \ref{as:g}. 

For the forthcoming optimization algorithm, Equation \eqref{eq:saddle} is used to define an
upper bound on the norm of $\hat{\mu}$, stated in the following lemma. It uses the Slater
point $\bar{x}$ from Assumption \ref{as:g}.\ref{as:gii} and any lower bound on $f$
over $X$, denoted $f_{lower}$, which exists under Assumption \ref{as:fX}. 

\begin{lemma} \label{lem:mubound}
For $(\hat{x}, \hat{\mu})$ a saddle point of $L$, 
\begin{equation}
\hat{\mu} \in \M := \left\{\mu \in \R^m_{+} : \|\mu\|_1 \leq  \frac{f(\bar{x}) - f_{lower}}{\min\limits_{j \in [m]} \left\{-g_j(\bar{x})\right\}}\right\}.
\end{equation} 
\end{lemma}
\emph{Proof:} See Section II-A in \cite{hale15c}. \hfill $\blacksquare$

Using the fact that a saddle point of $L$ provides a solution
to Problem \ref{prob:ensemble}, the remainder of the paper focuses on finding such saddle points. 
The algorithm used for this purpose includes an iterative Tikhonov regularization 
with an asymptotically vanishing stepsize, and was stated in \cite{bakushinskii74} for deterministic
variational inequalities and later in \cite{poljak78} for stochastic problems. 
To help describe the communications and computations in the cloud-based system, we
provide the general deterministic form of
this algorithm now, though in Section \ref{sec:main} stochasticity is introduced when the algorithm is made joint differentially private. 
The saddle-point finding algorithm uses the coupled update equations
\begin{subequations}\label{eq:alg1}
\begin{align} 
x(k+1) &= \Pi_{X}\left[x(k) - \gamma_k\left(L_{x}(k) + \alpha_kx(k)\right)\right] \label{eq:alg1p}\\
\mu(k+1) &= \Pi_{M}\left[\mu(k) + \gamma_k\left(L_{\mu}(k) - \alpha_k\mu(k)\right)\right], \label{eq:alg1d}
\end{align} 
\end{subequations}
where $\Pi_{X}$ and $\Pi_{M}$ are the Euclidean projections onto $X$ and $M$, respectively.
This update law will be referred to as Update \eqref{eq:alg1}. 
Here, $\gamma_k$ is a stepsize and $\alpha_k$ is the regularization parameter, and the values of both
will be provided in Theorem~\ref{thm:noisycon} in Section~\ref{sec:main}. 

\begin{figure}
\centering
\includegraphics[width=3.7in]{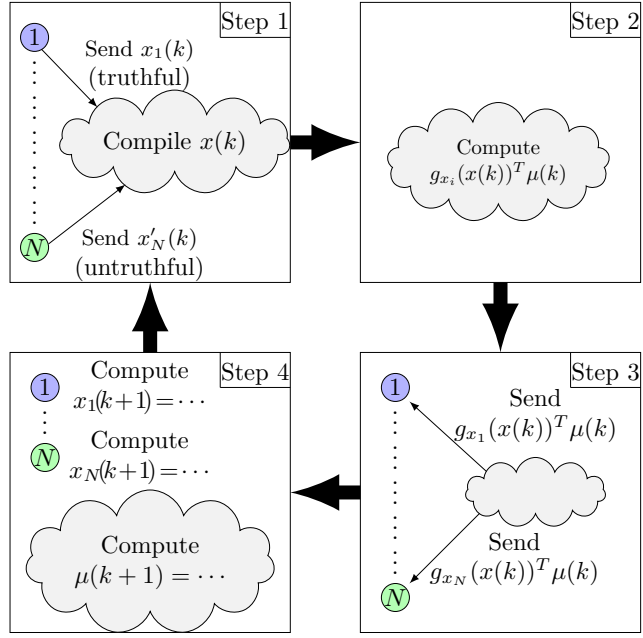}
\caption{The four steps of a communications cycle in the cloud-based system. First, each agent
sends its state to the cloud. 
Second, the cloud performs centralized computations required by the agents. 
Third, the cloud sends the results of these computations
to the agents. Fourth, agent $i$ computes $x_i(k+1)$ and the cloud computes $\mu(k+1)$ and 
this process repeats. In Step 1 depicted here, agent $N$ is misreporting its state to the cloud
by sending $x'_N(k)$ instead of $x_N(k)$. As the algorithm progresses 
this untruthful state propagates through the system. 
}
\label{fig:4steps}
\end{figure}

\subsection{Optimizing via the Cloud}
We now elaborate on the cloud-based architecture to be used and the means of executing Update~\eqref{eq:alg1} 
on the cloud and agents.  
To enforce joint differential privacy, the agents do not directly share any information with each other. 
Instead, the
agents route messages through a trusted cloud computer which aggregates all states in the network, performs computations
involving these states, and makes the results of these computations joint differentially private before sending them to the agents.

Splitting Equation \eqref{eq:alg1p} into each agent's states gives the 
update
\begin{equation} 
x_i(k+1) = \Pi_{X_i}\big[x_i(k) - \gamma_k\left(L_{x_i}(k) + \alpha_kx_i(k)\right)\big]
\end{equation}
for agent $i$, 
where expanding the term containing $L_{x_i}$ gives
\begin{equation} \label{eq:agenti}
x_i(k+1) = \Pi_{X_i}\Big[x_i(k) - \gamma_k\left(\nabla_{x_i} f_i\big(x_i(k)\big) + g_{x_i}\big(x(k)\big)^T\mu(k) + \alpha_kx_i(k)\right)\Big].
\end{equation}
Importantly, the right-hand side of Equation \eqref{eq:agenti} contains two terms which agent $i$ cannot compute on its own: $g_{x_i}\big(x(k)\big)$,
because it is a function of every state in the network, and $\mu(k)$, because its update law relies on $g\big(x(k-1)\big)$, which is also a function of every state
in the network. For this reason, the cloud computer is used to 
compute $\mu(k)$ and $g_{x_i}\big(x(k)\big)$ for each $i \in [N]$
at all timesteps $k$. 

Four actions occur within timestep $k$. 
First, agent $i$ sends $x_i(k)$ to the cloud and the cloud assembles
the ensemble state vector $x(k)$. Second, the cloud computes $g_{x_i}\big(x(k)\big)$ for every $i \in [N]$.
Third, the cloud sends $g_{x_i}\big(x(k)\big)^T\mu(k)$ to agent $i$. Fourth, the cloud
computes $\mu(k+1)$ and simultaneously agent $i$ computes $x_i(k+1)$, and then this process repeats.
This exchange of information is depicted in Figure \ref{fig:4steps}.
In the upper-left panel of Figure \ref{fig:4steps}, agents $1$ through $N-1$ report their
actual state to the cloud, while agent $N$ misreports its state to the cloud by sending some
$x'_N(k)$ instead of $x_N(k)$.
The arrows connecting boxes indicate how misreported states propagate through the system, eventually affecting
all agents' states.  
The cloud's computations in Steps 2 and 4 will be modified in Section~\ref{sec:jdp} to introduce
joint differential privacy, 
though the overall communications structure will remain the same. 

To reflect that agent $i$ receives the vector $g_{x_i}\big(x(k)\big)^T\mu(k)$ from the
cloud without knowing $g_{x_i}\big(x(k)\big)$ or $\mu(k)$ individually, agent $i$'s update
law is rewritten as
\begin{equation}
x_i(k+1) = \Pi_{X_i}\Big[x_i(k) - \gamma_k\left(\nabla_i f_i\big(x_i(k)\big) + q_i(k) + \alpha_kx_i(k)\right)\Big],
\end{equation}
where $\R^{n_i} \ni q_i(k) := g_{x_i}\big(x(k)\big)^T\mu(k)$. As in Equation~\eqref{eq:alg1d}, the cloud computes the dual update
according to
\begin{equation}
\mu(k+1) = \Pi_{M}\left[\mu(k) + \gamma_k\left(L_{\mu}(k) - \alpha_k\mu(k)\right)\right],
\end{equation}

\subsection{Full Problem Statement}
We now state the algorithm that will be used throughout the remainder of the paper. Below that, 
we identify the potential for misreporting states in this algorithm and state the problem
that is later solved using joint differential privacy. 

\begin{talgorithm} \label{alg:main}
$ $ \\
\emph{Step 0:} For all $i \in [N]$, initialize agent $i$ with $x_i(0)$, $X_i$, $f_i$, $\{\alpha_k\}_{k \in \N}$, and $\{\gamma_k\}_{k \in \N}$. 
Initialize the cloud with $\bar{x}$, $g$, $f_{lower}$, $\{\alpha_k\}_{k \in \N}$, and $\{\gamma_k\}_{k \in \N}$. Let the cloud compute
$M$ before the system begins optimizing. Set $k = 0$. \\
\emph{Step 1:} For all $i \in [N]$, the cloud computes $q_i(k)$ and sends it to agent $i$. \\
\emph{Step 2:} Agent $i$ computes 
\begin{equation}
x_i(k+1) = \Pi_{X_i}\Big[x_i(k) - \gamma_k\left(\nabla_i f_i\big(x_i(k)\big) + q_i(k) + \alpha_kx_i(k)\right)\Big], \tag{\value{equation}a}
\end{equation}
and sends a state value $x_i(k+1)$ to the cloud. \\
\emph{Step 3:} The cloud computes
\begin{equation}
\mu(k+1) = \Pi_{M}\left[\mu(k) + \gamma_k\left(L_{\mu}(k) - \alpha_k\mu(k)\right)\right]. \tag{\value{equation}b}
\end{equation}
\emph{Step 4:} Set $k := k+1$ and return to Step 1. \hfill $\Diamond$ 
\end{talgorithm}

\stepcounter{equation}

In Step~$2$ of Algorithm~\ref{alg:main}, agent $i$ may report a false state value to the cloud,
sending some $\tilde{x}_i(k+1) \neq x_i(k+1)$. This is the misreporting behavior that this paper seeks to prevent
using the cloud.  
The only influence the cloud has upon agent $i$
is through $q_i(k)$, and therefore the cloud must compute $q_i(k)$ in a way that incentivizes
agent $i$ to honestly report its state. 
The incentivization of this behavior is stated as Problem~\ref{prob:full} below, and the remainder
of the paper focuses on solving Problem~\ref{prob:full}. 

\begin{problem} \label{prob:full}
Execute Algorithm~\ref{alg:main} with the cloud computing $q_i(k)$ in a way that incentivizes
agent $i$ to honestly report $x_i(k+1)$ in Step~$2$, while still converging to a minimum. \hfill $\lozenge$
\end{problem}

It is assumed that each agent ultimately wants the constraints $g$ to be
satisfied, as would be the case when $g$ corresponds to some mission-critical
conditions that must be satisfied by the agents. 
However, an agent may wish to 
reduce the impact $g$ has upon its own state in order to reduce its cost. 
One way of affecting $g$ for this purpose is by reporting 
false states to the cloud over time; because all agents
want $g$ to be satisfied, a misreporting agent will still use $q_i(k)$
from the cloud in its state updates, but an agent can substantially
influence these messages for its own benefit through misreport.
In response to misreported states, other agents' messages from the cloud
in Update \eqref{eq:alg1} will be affected in a way that compensates
 for agent $i$'s false
reported states, thereby resulting in an unfair
distribution of the burden of $g$. 

This behavior cannot be
detected by the other agents or the cloud because only agent $i$ knows
$f_i$ and $X_i$, making no other entity in the network capable of
determining what agent $i$'s state should be (cf. Equation \eqref{eq:alg1p}).
Therefore, rather than detecting
manipulation of $g$, we seek to prevent
this behavior. Joint differential privacy provides a framework
for incentivizing truthful sharing of information and it is used here
to prevent the agents from manipulating $g$. 

Because agents are opportunistic but not malicious, they may send untruthful states
to the cloud, but they will not send states that harm the system or prevent
convergence of Algorithm~\ref{alg:main}. As a result, a misreported
state trajectory will have some relationship to an agent's true state trajectory and
this relationship is used in defining adjacency of signals in our 
joint differential privacy implementation. 
The next section details the manner in which noise is added to $q_i(k)$ using joint differential privacy, 
and Section \ref{sec:main} shows that this noise still allows for 
Algorithm~\ref{alg:main} to reach a minimum.

\section{Joint Differential Privacy} \label{sec:jdp}
This section recalls necessary details from 
both ordinary differential privacy and
joint differential privacy. Then it presents the joint differential privacy mechanism
that will be implemented on the preceding cloud architecture. 
It is critically important to note that while this section discusses privacy,
the ultimate goal is to apply joint differential privacy to induce truthful sharing
of states. All of the content on privacy in this section should therefore
be understood as making progress toward inducing truthful behavior.

\subsection{Differential Privacy Background}
Differential privacy as described in \cite{dwork06c} was originally designed to keep individual
database entries private whenever a database is queried, and
this is done by adding noise
to the responses to such queries. This idea was extended to dynamical
systems in \cite{leny14} in order to keep inputs to a system private
from anyone observing the outputs of that system. It is the dynamical
systems notion of differential privacy that is used below. 

The key idea behind differential privacy is that noise is added 
to make ``adjacent'' inputs produce ``similar'' outputs, and these notions
are made rigorous below. 
As above, let there be $N$ users, with the $i^{th}$ user contributing
a signal $u_i \in \uinp$. The space $\uinp$ is the space of sequences of $s_i$-vectors
in which every finite truncation of every element has finite $p_i$-norm. 
More explicitly, with $u_i(k) \in \R^{s_i}$ denoting
the $k^{th}$ element of the signal $u_i$, 
define the truncation operator $\trunc$ according to
\begin{equation} \label{eq:inpdef}
\trunc u = \begin{cases} u(k) & k \leq t \\ 0 & k > t \end{cases}.
\end{equation}
Then $u_i \in \uinp$ if and only if $\trunc u_i \in \lspi$
for all $t \in \N$. 
The full input space 
is then defined by the Cartesian product
$\inp = \prod_{i=1}^{n} \uinp$.
This paper focuses on the case where $p_i = 1$ for all $i \in [N]$. 

To formalize the notion of adjacency of inputs in $\inp$, we define a binary, 
symmetric adjacency relation $\adj_{B}^i$, which is parameterized by $B > 0$ and 
an index $i \in [N]$.
Its definition uses the notation
$u_{-i} = (u_1, \ldots, u_{i-1}, u_{i+1}, \ldots, u_N)$. 
The adjacency relation takes the form
$\adj^i_{B} : \inp \times \inp \to \{0, 1\}$ and has the following definition,
as stated in \cite{leny14}. 
\begin{definition} \label{def:adjacency}
Two inputs $u$, $\tilde{u} \in \inp$ satisfy ${\adj^i_{B}(u, \tilde{u}) = 1}$ 
for some $i \in [N]$ if and only if
$\|u_i - \tilde{u}_i\|_{p_i} \leq B \textnormal{ and } u_{-i} = \tilde{u}_{-i}$.
The signals $u$ and $\tilde{u}$ are then said to be adjacent with respect to agent $i$. 
If $i$ is arbitrary, the relation $\adj_{B}$ is used, and $u$ and $\tilde{u}$ are simply called adjacent. 
\hfill $\blacklozenge$
\end{definition}

This section maintains use of the symbol $u$ for system inputs (rather than $x$ as will be in subsequent sections) to maintain
continuity with the references cited here for private dynamical systems. 

Inputs from $\inp$ are assumed to pass into a causal, deterministic system $\sys$, which produces
outputs in $\outp$. 
To define when two outputs are ``similar,'' the notion of a \emph{mechanism} is used. 
In the context of private dynamical systems, 
a mechanism is a means of adding noise to an otherwise deterministic system
to provide privacy to that system's inputs. 
Formally, for a fixed a probability space $(\Omega, \mathcal{F}, \p)$,
a mechanism is a map of the form
\begin{equation}
M : \inp \times \Omega \to \outp.
\end{equation}
The mechanism $M$ must provide differential privacy to its input trajectories, which
fundamentally means that the output of $M$ should be insensitive to changes in
its inputs. Differential privacy captures this notion by requiring that,
whenever $\adj_{B}(u, \tilde{u})$ holds,
the probability distributions of $Mu$ and $M\tilde{u}$ satisfy
\begin{equation}
\p(Mu \in S) \leq e^{\epsilon}\p(M\tilde{u} \in S)
\end{equation}
for all $S$ in an appropriate $\sigma$-algebra. 

Using a result from the literature, we now state a finite-time criterion
which holds if and only if $M$ keeps entire trajectories private.
Below, the notation
$v_{0:k} := \big(v(0), \ldots, v(k)\big)$
is used to refer to the first $k+1$ entries of $v \in \outp$. 
\begin{lemma} \label{lem:dp}
Let $\epsilon \geq 0$ be given. 
For a dynamical system, a mechanism $M$ is $\epsilon$-differentially private if and only if, for all $u$, $\tilde{u}$ satisfying $\adj_{B}(u, \tilde{u}) = 1$ and
for all times $k$,
\begin{equation} \label{eq:dpdef}
\p\big((M(u))_{0:k} \in A\big) \leq e^{\epsilon}\p\big((M(\tilde{u}))_{0:k} \in A\big) \textnormal{ for all } A \in \mathscr{B}^{(k+1)r},
\end{equation}
where $\mathscr{B}^{d}$ is the Borel sigma-algebra on $\R^d$ and $r$ is the dimension of the output space. 
\end{lemma} 
\emph{Proof:} See Lemma 2 in \cite{leny14}. \hfill $\blacksquare$

In Lemma \ref{lem:dp}, the value of $\epsilon$ determines the level of privacy
afforded to the input signals, and decreasing its value leads to improved privacy at
the cost of adding higher variance noise. 
Typical values of $\epsilon$ in the literature range from $0.1$ to $\ln 3$. 

\subsection{Joint Differential Privacy}
We now elaborate on the application of joint differential privacy to Problem~\ref{prob:full}. 
To promote truth-telling behaviors, limits are imposed on 
the ability of any agent to reduce its cost by reporting a false state trajectory to the cloud.
These limits are enforced using 
joint differential privacy, which is a relaxation of ordinary differential privacy for use in multi-agent systems,
and it will be shown that this framework is sufficient for the goal of reducing an agent's ability
to benefit from misreporting its state. 
For joint differential privacy, the ``system'' of interest is comprised by the computations the cloud
carries out in accordance with
Update~\eqref{eq:alg1}, and this point is discussed further below. 
For now it suffices to point out that the output of this system is a tuple of private forms
of all $q_i$'s. Denoting the private form of $q_i(k)$ by $\tilde{q}_i(k)$ (whose exact
form will be given later), the output of the cloud at time $k$ is denoted by
\begin{equation}
y(k) = \big(\tilde{q}_1(k), \ldots, \tilde{q}_N(k)\big). 
\end{equation}

Let $\Mjdp$ denote a mechanism for joint
differential privacy and let $\Mjdp_{-i}$ denote the same mechanism with the $i^{th}$ output removed, i.e.,
the output of $\Mjdp_{-i}$ is
\begin{equation}
y_{-i}(k) = \big(\tilde{q}_1(k), \ldots, \tilde{q}_{i-1}(k), \tilde{q}_{i+1}(k), \ldots, \tilde{q}_N(k)\big). 
\end{equation}
A mechanism $\Mjdp$ is joint differentially private if, for any $i \in [N]$, $\Mjdp_{-i}$ preserves differential
privacy for inputs adjacent with respect to $i$. Joint differential privacy
for databases has been defined in \cite{kearns12}, though, to our knowledge,
it has not yet been used for dynamical systems. 
Using Lemma 2 in \cite{leny14}, the following lemma states a
finite-time criterion for trajectory-level
joint differential privacy. 
\begin{lemma} (Joint differential privacy for dynamical systems) \label{lem:jdp}
Let the privacy parameter $\epsilon \geq 0$ be given and let 
$\Mjdp$ be a mechanism whose output is an $N$-tuple. Then $\Mjdp$ is $\epsilon$-joint differentially
private if and only if, for any $i \in [N]$, for all $u, \tilde{u} \in \inp$ satisfying
$\adj^i_{B}(u, \tilde{u}) = 1$, all times $k$, and all $A \in \mathscr{B}^{(k+1)(n - n_i)}$, $\Mjdp$ satisfies
\begin{equation}
\p\big((\Mjdp_{-i}(u))_{0:k} \in A\big) \leq e^{\epsilon}\p\big((\Mjdp_{-i}(\tilde{u}))_{0:k} \in A\big) 
\end{equation}
where $\mathscr{B}^d$ is the Borel sigma-algebra on $\R^d$. \hfill $\blacksquare$
\end{lemma}

Joint differential privacy guarantees that when agent $i$'s input changes by 
a small amount, the outputs corresponding to other agents do not change by much. 
A useful characteristic of both ordinary and joint differential privacy is 
their resilience to post-processing, which guarantees that post-hoc transformations of private
data cannot weaken the privacy guarantees afforded to that data. 
This result is formalized below. 

\begin{lemma} (Resilience to post-processing; \cite[Proposition 2.1]{dwork13}) \label{lem:resilience}
Let $M$ be an $\epsilon$-differentially private mechanism and let $f$ be a function such that the
composition $f \circ M$ is well-defined. Them $f \circ M$ is also $\epsilon$-differentially private. \hfill $\blacksquare$
\end{lemma}

\subsection{The Laplace Mechanism}
To enforce differential privacy for a particular choice of $\epsilon$, 
noise must be added somewhere in the system, and the distribution of
this noise must be determined. 
One common mechanism in the literature
draws noise from the Laplace distribution, used in both \cite{dwork06c} and \cite{leny14}, and
this mechanism is used 
to provide $\epsilon$-joint differential privacy.
To define this mechanism, the notion of the $\ell_p$ sensitivity of a system is now introduced. 
\begin{definition} \label{def:sensivitiy}
The $\ell_p$ sensitivity of a deterministic, causal system $\sys$ is defined as
\begin{equation}
\Delta_p\sys = \sup_{u, \tilde{u} : \adj_{B}(u, \tilde{u}) = 1} \left\|\sys(u) - \sys(\tilde{u})\right\|_{\ell_p}. \tag*{$\blacklozenge$}
\end{equation} 
\end{definition}

The Laplace mechanism is stated in terms of the $\ell_1$ sensitivity of a system. Below, 
the notation $\lap(b)$ denotes a scalar Laplace distribution with mean zero and scale parameter
$b$, i.e.,
\begin{equation}
\lap(b) := \frac{1}{2b}\exp\left(-\frac{|x|}{b}\right). 
\end{equation} 
\begin{lemma} (Laplace mechanism; \cite[Theorem 4]{leny14}) \label{lem:laplace}
Let $\epsilon \geq 0$ be given. The Laplace mechanism defined by
\begin{equation}
M(u) = \sys(u) + w
\end{equation}
with $w(k) \sim \lap(b)^r$ is $\epsilon$-differentially private for $b \geq \Delta_1\sys/\epsilon$, with
$r$ the dimension of the output space of the system. \hfill $\blacksquare$
\end{lemma}

Lemma \ref{lem:laplace} says that one can make a system private by adding noise drawn from a Laplace
distribution to that system's output at each timestep. This idea is implemented for Problem~\ref{prob:full}
in the next section.

\section{Optimizing Under Joint Differential Privacy} \label{sec:main}
In this section, the privacy results in Section~\ref{sec:jdp} are applied 
to Problem~\ref{prob:full}, and Algorithm~\ref{alg:main} is made 
joint differentially private. 
Formally, $g$ and $g_{x_i}$ are treated as memoryless
dynamical systems, and 
joint differential privacy is used to ensure that
they keep the agents' state trajectories, which are the inputs to these systems, private. 

\subsection{Stochastic Optimization Algorithm}
From Lemma \ref{lem:jdp}, we see that enforcing joint differential privacy
for the states in the network requires that the cloud make
$q_i(k):=g_{x_i}\big(x(k)\big)^T\mu(k)$ private before it is sent to agent $i$. 
To make $g_{x_i}(x(k))$ private, noise can be added to it directly
and this is done below. To make $\mu(k)$ private, we use the fact
that computing $\mu(k)$
relies on $g(x(k-1))$, and add noise to make $g(x(k-1))$ private. Then computing $\mu(k)$ is joint differentially
private by the post-processing property in Lemma \ref{lem:resilience}. 
Similarly, if the noisy forms of both $g_{x_i}(x(k))$ and $\mu(k)$ are joint differentially private,
their product is as well, again by Lemma~\ref{lem:resilience}. 
Adding noise in this way, Algorithm~\ref{alg:main} is modified to state the
joint differentially private optimization algorithm below.
\begin{talgorithm} \label{alg:noisy}
$ $ \\
\emph{Step 0:} For all $i \in [N]$, initialize agent $i$ with $x_i(0)$, $X_i$, $f_i$, $\{\alpha_k\}_{k \in \N}$, and $\{\gamma_k\}_{k \in \N}$. 
Initialize the cloud with $\bar{x}$, $g$, $f_{lower}$, $\{\alpha_k\}_{k \in \N}$, and $\{\gamma_k\}_{k \in \N}$. Let the cloud compute
$M$ before the system begins optimizing. Set $k = 0$. \\
\emph{Step 1:} For all $i \in [N]$, the cloud computes \\ ${\tilde{q}_i(k) := \big(g_{x_i}\big(x(k)\big) + w_i(k)\big)^T\mu(k)}$ and sends it to agent $i$. \\
\emph{Step 2:} Agent $i$ computes 
\begin{equation}
x_i(k+1) = \Pi_{X_i}\Big[x_i(k) - \gamma_k\left(\nabla_i f_i\big(x_i(k)\big) + \tilde{q}_i(k) + \alpha_kx_i(k)\right)\Big], \tag{\value{equation}a}
\end{equation}
and sends a state value $x_i(k+1)$ to the cloud. \\
\emph{Step 3:} The cloud computes
\begin{equation}
\mu(k+1) = \Pi_{M}\left[\mu(k) + \gamma_k\left(L_{\mu}(k) + w_g(k) - \alpha_k\mu(k)\right)\right]. \tag{\value{equation}b}
\end{equation}
\emph{Step 4:} Set $k := k+1$ and return to Step 1. \hfill $\Diamond$ 
\end{talgorithm}

To solve Problem~\ref{prob:full}, 
Algorithm~\ref{alg:noisy} must implement joint differential privacy using the noise terms $w_g$ and $w_i$, while
still converging to a minimum. The following theorem gives conditions on $w_g$ and each $w_i$ under
which convergence to a minimum is guaranteed. Then Theorem~\ref{thm:mech} shows that joint
differential privacy is achieved under these conditions, and Section~\ref{sec:eta} shows that
each agent's incentive for misreport is indeed limited due to joint differential privacy. 

\begin{theorem} \label{thm:noisycon}
Let $(\hat{x}, \hat{\mu})$ denote the least-norm saddle point of $L$. 
Algorithm \ref{alg:noisy} satisfies
\begin{equation}
\lim_{k \to \infty} \mathbb{E}\big[\|x(k) - \hat{x}\|_2^2\big] = 0 \textnormal{ and } \lim_{k \to \infty} \mathbb{E}\big[\|\mu(k) - \hat{\mu}\|_2^2\big] = 0
\end{equation}
if i. $w_i(k)$ and $w_g(k)$ have zero mean for all $k$ \\ 
ii. $\alpha_k = \bar{\alpha}k^{-c_1}$ and $\gamma_k = \bar{\gamma}k^{-c_2}$, where $0 < c_1 < c_2$, $c_1 + c_2 < 1$, $0 < \bar{\alpha}$, and $0 < \bar{\gamma}$.
\end{theorem}
\emph{Proof:} See \cite[Theorem 6]{poljak78}. \hfill $\blacksquare$

It remains to be shown that Condition i of Theorem~\ref{thm:noisycon} can be satisfied when joint differential
privacy is implemented, and this is done next. 

\subsection{Calibrating Noise for Joint Differential Privacy}
Here the systems being made private are $g$ and $g_{x_i}$, and 
the mechanisms acting for joint differential privacy add noise to $g(x(k-1))$ when computing
$\mu(k)$ and add noise to $g_{x_i}(x(k))$ when computing $\tilde{q}_i(k)$. 
It was shown in Section~\ref{sec:jdp} that the noise added in Algorithm \ref{alg:noisy} 
will enforce $\epsilon$-joint differential privacy as long as it has large enough variance. 
To determine the variance of noise that must be added by the Laplace mechanism, bounds are derived on the $\ell_1$ sensitivity of each $g_{x_i}$ and $g$ below. 
\begin{lemma} \label{lem:sensitivity}
For the relation $\adj_{B}$, the $\ell_1$ sensitivities of $g_{x_i}$ and $g$ satisfy
$\Delta_1 g_{x_i} \leq \lipgradgi B$ and $\Delta_1 g \leq \lipg B$.
\end{lemma} 
\emph{Proof:} See \cite{hale15c}. 
\hfill $\blacksquare$

Using Lemmas \ref{lem:laplace} and \ref{lem:sensitivity}, we see that if $w_i(k) \sim \Lap(b_i)$ 
with $b_i \geq \Delta_1g_{x_i}/\epsilon$ and
$w_g(k) \sim \Lap(b_g)$ with $b_g \geq \Delta_1 g/\epsilon$ for all $k$, then all states are afforded (ordinary) $\epsilon$-differential privacy
in Algorithm \ref{alg:noisy}. It turns out that this privacy and its resilience
to post-processing imply that $\epsilon$-joint differential privacy holds as well, which is stated formally in the following theorem. 
\begin{theorem} \label{thm:mech}
Consider the mechanism $\Mjdp$ defined by
\begin{equation}
\Mjdp\big(x(k)\big) = (\tilde{q}_1(k), \ldots, \tilde{q}_N(k)),
\end{equation}
with $\tilde{q}_i(k) := \big(g_{x_i}\big(x(k)\big) + w_i(k)\big)^T\mu(k)$
as defined in Algorithm~\ref{alg:noisy} and $\mu(k)$ computed
as in Algorithm~\ref{alg:noisy}. 
At each time $k$, if $w_i(k) \sim \Lap(b_i)$ 
with $b_i \geq \Delta_1 g_{x_i}/\epsilon$ and
$w_g(k) \sim \Lap(b_g)$ with $b_g \geq \Delta_1 g/\epsilon$,
then $\Mjdp$ is $\epsilon$-joint differentially
private.
\end{theorem}
\emph{Proof:} We examine $\Mjdp_{-i}$ for an arbitrary $i \in [N]$ whose output at time $k$ is
\begin{equation}
\Mjdp_{-i}\big(x(k)\big) = \left(\tilde{q}_1(k), \ldots, \tilde{q}_{i-1}(k), \tilde{q}_{i+1}(k), \ldots, \tilde{q}_N(k)\right).
\end{equation}
Examining some $j \in [N] \backslash \{i\}$, the $j^{th}$ output $\Mjdp_j\big(x(k)\big)$ is
\begin{equation}
\Mjdp_j\big(x(k)\big) = \tilde{q}_j(k) := \left(g_{x_j}\big(x(k)\big) + w_j(k)\right)^T\mu(k).
\end{equation}
In light of the fact that
$w_g(k-1) \sim \Lap(b_g)$ with 
$b_g \geq \Delta_1 g/\epsilon$, 
we see that
$\mu(k)$ keeps $\bx(k-1)$ $\epsilon$-differentially private. 
Similarly, 
$g_{x_j}\big(x(k)\big) + w_j(k)$ keeps $x(k)$ $\epsilon$-differentially private because $w_j(k) \sim \lap(b_j)$ and
$b_j \geq \Delta_1 g_{x_j}/\epsilon$.
By Lemma \ref{lem:resilience}, the product $\tilde{q}_j(k) = \left(g_{x_j}\big(x(k)\big) + w_j(k)\right)^T\mu(k)$
keeps $\bx$ $\epsilon$-differentially private because it is the result of post-processing
two differentially private quantities. 
By the same reasoning,
\begin{equation}
\Mjdp_{-i}\big(x(k)\big) = \big(\tilde{q}_1(k), \ldots, \tilde{q}_{i-1}(k), \tilde{q}_{i+1}(k), \ldots, \tilde{q}_N\big)
\end{equation}
simply post-processes differentially private information. 
From Lemma \ref{lem:jdp}
we conclude that $\Mjdp$ is $\epsilon$-joint differentially private. 
\hfill $\blacksquare$

The next section describes the application
of this mechanism to inducing approximately-truthful behavior 
in multi-agent optimization through the computation of $\beta$-approximate minima.

\section{Computing $\beta$-approximate Minima} \label{sec:eta}
This section presents the main result of the paper: joint differential
privacy results in there being
only minimal incentive for an agent to misreport its state to the cloud. 
Toward showing this result, 
a uniform upper bound on $f_i$ over $X_i$ is first presented.  
\begin{lemma} \label{lem:fmax}
Let $\bar{x}$ denote a Slater point for $g$. 
Then for all $i \in [N]$,
$f_i(x_i) \leq \lambda_i$
for all $x_i \in X_i$, where $\lambda_i := f_i(\bar{x}_i) + \lipfi D_i$.
\end{lemma}
\emph{Proof:} Using the Mean-Value Theorem and the Cauchy-Schwarz inequality, we have
\begin{align}
f_i(x_i) &= f_i(\bar{x}_i) + \nabla f_i(z_i)^T(x_i - \bar{x}_i) \\
         &\leq f_i(\bar{x}_i) + \|\nabla_i f_i(z_i)\|\cdot\|x_i - \bar{x}_i\|,
\end{align}
for some $z_i \in X_i$. The result follows by bounding $\|\nabla_i f_i(z_i)\|$ by $\lipfi$ and
 bounding $\|x_i - \bar{x}_i\|$ by $D_i$. \hfill $\blacksquare$

The next lemma bounds the difference in $f_i$ at points along
two feasible state trajectories for agent $i$. 
\begin{lemma} \label{lem:costdiff}
For any $\sone, \stwo \in \uinp$ and any time $k$, one finds
\begin{equation}
\left|f_i\big(\sone(k)\big) - f_i\big(\stwo(k)\big)\right| \leq \fgaini := \min\{\lipfi D_i, 2\lambda_i\}. 
\end{equation}
\end{lemma}
\emph{Proof:} Using the Lipschitz property of $f_i$,
\begin{equation}
\left|f\big(\sone(k)\big) - f\big(\stwo(k)\big)\right| \leq \lipfi \|\sone(k) - \stwo(k)\|_1 \leq \lipfi D_i.
\end{equation}
On the other hand, the triangle inequality gives
\begin{equation}
\left|f\big(\sone(k)\big) - f\big(\stwo(k)\big)\right| \leq |f\big(\sone(k)\big)| + |f\big(\stwo(k)\big)| \leq 2\lambda_i,
\end{equation}
where the second inequality follows from Lemma \ref{lem:fmax}. 
\hfill $\blacksquare$

The main result of the paper is now presented. 
Below, each expected value
is over the randomness introduced by the mechanism $\Mjdp$.
For clarity, this theorem tracks the state agent $i$ has reported
to the cloud and the state trajectory of every
\emph{other} agent. 
The notation
$\E[f_i\big(x_i(k)\big) | y_i, v_{-i}]$
is used to denote agent $i$'s cost at time $k$ when agent $i$ has reported the trajectory $y_i$ to the cloud
and every other agent has followed the trajectory $v$.
The symbol $x_i$ is always used as the argument to $f_i$ because
$f_i$ always depends on the true state of agent $i$, not the state
it reports. 

\begin{theorem} \label{thm:eta} 
Suppose Assumptions \ref{as:fX} and \ref{as:g} hold. 
Let the agents and cloud execute Algorithm \ref{alg:noisy} 
with the cloud implementing the mechanism $\Mjdp$. 
Then for $\epsilon \in (0, 1)$, 
all agents sharing their true states in
Algorithm \ref{alg:noisy} results in a $\beta$-approximate minimum. 
In particular, at all times $k$ and 
for any state trajectories $\bx,\btx \in \inp$ 
satisfying $\adj_{B}^i(\bx, \btx) = 1$,
we have
\begin{equation}
\E[f_i(x_i(k)) | x_i, x_{-i}] \leq \E[f_i(x_i(k)) | x_i', \tilde{x}_{-i}] + \beta,
\end{equation}
where $\beta = 2\max_{i \in [N]} \fgaini + \epsilon\lambda_i$, and where
$x_i'$ is a misreported state trajectory from agent $i$. 
\end{theorem}
\emph{Proof:} From Lemma \ref{lem:costdiff} we have
\begin{equation} \label{eq:main1}
\E[f_i(x_i(k)) | x_i, x_{-i}] \leq \E[f_i(x_i(k)) | \tilde{x}_i, x_{-i}] + \fgaini.
\end{equation} 
Using Theorem \ref{thm:mech} and the definition of 
$\epsilon$-joint differential privacy we find
\begin{equation}
\E[f_i(x_i(k)) | \tx_i, x_{-i}] \leq e^{\epsilon}\E[f_i(x_i(k)) | \tx_i, \tx_{-i}], 
\end{equation}
which we substitute into Equation \eqref{eq:main1} to find
\begin{equation}
\E[f_i(x_i(k)) | x_i, x_{-i}] \leq e^{\epsilon}\E[f_i(x_i(k)) | \tx_i, \tx_{-i}] + \fgaini.
\end{equation}
Using $e^{\epsilon} \leq 1 + 2\epsilon$ for $\epsilon \in (0,1)$ gives
\begin{equation} \label{eq:almostthere}
\E[f_i(x_i(\runtime)) | x_i, x_{-i}] \leq \E[f_i(x_i(\runtime)) | \tx_i, \tx_{-i}] + 2\epsilon\E[f_i(x_i(\runtime)) | \tx_i, \tx_{-i}] + \fgaini,
\end{equation}
were we apply Lemma \ref{lem:fmax} to get
\begin{equation} \label{eq:almostthere}
\E[f_i(x_i(\runtime)) | x_i, x_{-i}] \leq \E[f_i(x_i(\runtime)) | \tx_i, \tx_{-i}] + 2\epsilon\lambda_i + \fgaini.
\end{equation}
A second application of Lemma \ref{lem:costdiff} gives
\begin{equation}
\E[f_i(x_i(\runtime)) | \tx_i, \tx_{-i}] \leq \E[f_i(x_i(\runtime)) | x'_i, \tx_{-i}] + \rho_i,
\end{equation}
and substituting this inequality into Equation \eqref{eq:almostthere} gives 
\begin{equation}
\E[f_i(x_i(\runtime)) | x_i, x_{-i}] \leq \E[f_i(x_i(\runtime)) | x'_i, \tx_{-i}] + \beta,
\end{equation}
as desired. 
\hfill $\blacksquare$

While the $\rho_i$ term in $\beta$ is a feature of the problem itself, the $\epsilon \lambda_i$ term results
directly from the untruthfulness of agent $i$, and it is precisely this term which can be influenced
using the privacy parameter $\epsilon$, allowing a network operator to directly counteract the influence of false
information. Of course, shrinking $\epsilon$ requires that more noise be added which, in general, degrades
performance in the system. One must therefore balance the two objectives of incentivizing truthful information
sharing and system performance based upon the needs in a particular application.

\section{Simulation Results} \label{sec:simulation}
A simulation was run consisting of $N = 8$ agents
each with $x_i \in \R^2$, and $m = 4$ constraints.
The constraint $X_i = [-10, 10]^2$ for all $i \in [N]$, and
$f_i(x_i) = \frac{1}{2}\|x_i - t_i\|_2^2$,
for all $i \in [N]$; 
the value of each $t_i$ can be found in Table \ref{tab:constants}. 
The constraints used were 
\begin{equation}
g(x) = \left(\begin{array}{c} \|x_1 - x_2\|_2^2 + \|x_1 - x_3\|_2^2 - 5 \\
                              \|x_4 - x_5\|_2^2 + \|x_4 - x_6\|_2^2 - 3 \\
                              \|x_7 - x_8\|_2^2 + \|x_7 - x_6\|_2^2 - 3 \\
                              \|x_5 - x_3\|_2^2 + \|x_5 - x_7\|_2^2 - 5 \end{array}\right).         
\end{equation}

The value $B = 3$ was used for adjacency and the privacy parameter was chosen to be $\epsilon = \ln 3$.
The distributions of noise added are shown in Table \ref{tab:constants}; in
addition to those values, $w_g \sim \lap(327.69)$. The stepsize and
regularization parameters were chosen to be
$\gamma_k = 0.01k^{-3/5}$ and $\alpha_k = 0.5k^{-1/3}$. 
Two adjacent problems were run for $250,000$ timesteps each, with agent $6$ being
untruthful in one of them. In that run, agent $6$ reported its unconstrained minimizer $t_6$
to the cloud at each timestep instead of its actual state.  

%


\begin{figure}
\centering
\includegraphics[width=3.4in]{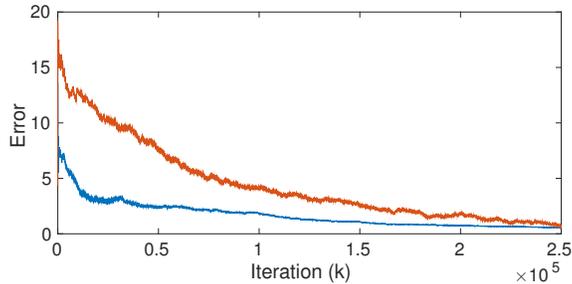}
\caption{A plot of the distance to the saddle point in Problem \ref{prob:ensemble}
in the primal space (lower curve) and dual space (upper curve) when all agents are truthful.}
\label{fig:error}
\end{figure}

\begin{figure}
\centering
\includegraphics[width=3.4in]{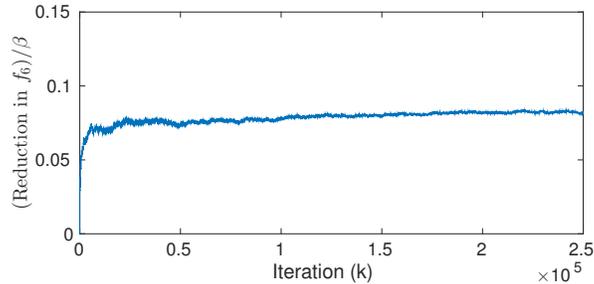}
\caption{A plot of the decrease in cost agent $6$ attains through misreporting
its state. The ordinate is normalized by $\beta$.}
\label{fig:eta}
\end{figure}

\begin{table}
\centering
\begin{tabular}{|c|c|c|} \hline
$i$ & $t_i$ & Distribution of $w_i$ \\ \hline\hline
$1$ & $[6 \,\,\, {-4}]^T$ & $\lap(10.92)$ \\ \hline
$2$ & $[2 \,\,\, 2]^T$ & $\lap(5.46)$ \\ \hline
$3$ & $[{-7} \,\,\, 7]^T$ & $\lap(5.46)$ \\ \hline
$4$ & $[8 \,\,\, {-9}]^T$ & $\lap(10.92)$ \\ \hline
$5$ & $[3 \,\,\, {-7}]^T$ & $\lap(16.38)$ \\ \hline
$6$ & $[10 \,\,\, 10]^T$ & $\lap(10.92)$ \\ \hline
$7$ & $[{-10} \,\,\, {-10}]^T$ & $\lap(16.38)$ \\ \hline
$8$ & $[6 \,\,\, {-6}]^T$ & $\lap(5.46)$ \\ \hline
\end{tabular}
\caption{Values of $t_i$ and distributions of $w_i$ for $i \in [8]$.}
\label{tab:constants}
\end{table}

Figure \ref{fig:error} plots the distance to the saddle point in the
primal and dual spaces when all agents are truthful. The final error values
were $\|x(250,000) - \hat{x}\|_2 = 0.5367$ and $\|\mu(250,000) - \hat{\mu}\|_2 = 0.6870$,
indicating close convergence in each space despite the large amount of noise in the system.

Figure \ref{fig:eta} plots the decrease in cost agent $6$ sees
in misreporting its state. This decrease
is not only bounded by $\beta$, but is bounded
above by $0.1\beta$ for all time, indicating that $\beta$ may be loose in
some cases; this is not surprising given that $\beta$ uses Lipschitz
constants and set diameters, which are ``worst-case'' in the sense that they
give maximum values over all possible states. Nonetheless, the algorithm
can be seen both to converge and compute a $\beta$-approximate minimum, indicating
that Problem~\ref{prob:full} has been solved.

\section{Conclusion}
It was shown that joint differential privacy can be used in multi-agent optimization
to incentivize truthful information sharing.
Applications of the work presented here include 
any multi-agent setting in which the iterates of an optimization algorithm
correspond to some physical quantity of interest. 
Future directions include allowing asynchronous
communications in order to account for systems with communication latency and poor
channel quality. 

\bibliographystyle{plain}{}
\bibliography{sources}

\begin{thebibliography}{10}

\bibitem{bakushinskii74}
AB~Bakushinskii and BT~Polyak.
\newblock Solution of variational inequalities.
\newblock {\em Doklady Akademii SSSR}, 219(5):1038--1041, 1974.

\bibitem{bertsekas03}
Dimitri~P Bertsekas, Angelia Nedic, and Asuman~E Ozdaglar.
\newblock {\em Convex analysis and optimization}.
\newblock Athena Scientific, 2003.

\bibitem{boyd11}
S.~Boyd, N.~Parikh, E.~Chu, B.~Peleato, and J.~Eckstein.
\newblock Distributed optimization and statistical learning via the alternating
  direction method of multipliers.
\newblock {\em Found. Trends Mach. Learn.}, 3(1), January 2011.

\bibitem{braynov04}
S.~Braynov and M.~Jadliwala.
\newblock Detecting malicious groups of agents.
\newblock In {\em Multi-Agent Security and Survivability, 2004 IEEE First
  Symposium on}, pages 90--99, Aug 2004.

\bibitem{cortes02}
Jorge Cortes, Sonia Martinez, Timur Karatas, and Francesco Bullo.
\newblock Coverage control for mobile sensing networks.
\newblock In {\em Robotics and Automation, 2002. Proceedings. ICRA'02. IEEE
  International Conference on}, volume~2, pages 1327--1332. IEEE, 2002.

\bibitem{dwork13}
C.~Dwork and A.~Roth.
\newblock The algorithmic foundations of differential privacy.
\newblock {\em Theoretical Computer Science}, 9(3-4):211--407, 2013.

\bibitem{dwork06c}
Cynthia Dwork, Frank McSherry, Kobbi Nissim, and Adam Smith.
\newblock Calibrating noise to sensitivity in private data analysis.
\newblock In {\em Proc. of the Third Conference on Theory of Cryptography},
  TCC'06, pages 265--284, Berlin, Heidelberg, 2006. Springer-Verlag.

\bibitem{fagiolini08}
A.~Fagiolini, M.~Pellinacci, G.~Valenti, G.~Dini, and A.~Bicchi.
\newblock Consensus-based distributed intrusion detection for multi-robot
  systems.
\newblock In {\em Robotics and Automation, 2008. ICRA 2008. IEEE International
  Conference on}, pages 120--127, May 2008.

\bibitem{hale15c}
M.T. Hale and M.~Egerstedt.
\newblock Cloud-enabled multi-agent optimization with constraints and
  differentially private states.
\newblock 2016.
\newblock Submitted for publication. Available at
  http://arxiv.org/abs/1507.04371.

\bibitem{han16}
S.~Han, U.~Topcu, and G.~J. Pappas.
\newblock Differentially private distributed constrained optimization.
\newblock {\em IEEE Transactions on Automatic Control}, PP(99):1--1, 2016.

\bibitem{han15}
Shuo Han, U.~Topcu, and G.J. Pappas.
\newblock An approximately truthful mechanism for electric vehicle charging via
  joint differential privacy.
\newblock In {\em American Control Conference (ACC), 2015}, July 2015.

\bibitem{hsu16}
Justin Hsu, Zhiyi Huang, Aaron Roth, and Zhiwei~Steven Wu.
\newblock Jointly private convex programming.
\newblock In {\em Proc. of the 27th Annual ACM-SIAM Symposium on Discrete
  Algorithms}, SODA '16, 2016.

\bibitem{kearns12}
M.~Kearns, M.~Pai, A.~Roth, and J.~Ullman.
\newblock Mechanism design in large games: Incentives and privacy.
\newblock {\em CoRR}, abs/1207.4084, 2012.

\bibitem{kuhn51}
H.~W. Kuhn and A.~W. Tucker.
\newblock Nonlinear programming.
\newblock In {\em Proc. of the 2nd Berkeley Symposium on Math. Stat. and
  Prob.}, Berkeley, Calif., 1951. University of California Press.

\bibitem{leny14}
J.~Le~Ny and G.J. Pappas.
\newblock Differentially private filtering.
\newblock {\em Automatic Control, IEEE Transactions on}, 59(2):341--354, Feb
  2014.

\bibitem{mcsherry07}
Frank McSherry and Kunal Talwar.
\newblock Mechanism design via differential privacy.
\newblock In {\em Proceedings of the 48th Annual IEEE Symposium on Foundations
  of Computer Science}, FOCS '07, pages 94--103, Washington, DC, USA, 2007.
  IEEE Computer Society.

\bibitem{nazari14}
M.H. Nazari, Z.~Costello, M.J. Feizollahi, S.~Grijalva, and M.~Egerstedt.
\newblock Distributed frequency control of prosumer-based electric energy
  systems.
\newblock {\em Power Systems, IEEE Transactions on}, 29(6), Nov 2014.

\bibitem{poljak78}
BT~Poljak.
\newblock Nonlinear programming methods in the presence of noise.
\newblock {\em Math. programming}, 14(1):87--97, 1978.

\bibitem{soltero13}
Daniel~E Soltero, Mac Schwager, and Daniela Rus.
\newblock Decentralized path planning for coverage tasks using gradient descent
  adaptive control.
\newblock {\em The International Journal of Robotics Research}, 2013.

\end{thebibliography}

\end{document}